\documentclass[12pt]{article}

\usepackage{amsfonts}
\usepackage{latexsym}
\usepackage{amsmath}
\usepackage{amssymb}

\newcommand{\veps}{\varepsilon}

\newcommand{\Q}{\mathbb{Q}}
\newcommand{\card}{\mathrm{card}}
\newcommand{\C}{\mathbb{C}}
\newcommand{\N}{\mathbb{N}}

\newcommand{\abs}[1]{\left|#1\right|}
\newcommand{\set}[1]{\left\{#1\right\}}
\newcommand{\brac}[1]{\left(#1\right)}
\newcommand{\sqbrac}[1]{\left[#1\right]}
\newcommand{\ord}{\textrm{ord}}

\newcommand{\di}{\textrm{di}}

\newtheorem{lettertheorem}{Theorem}
\newtheorem{letterlemma}[lettertheorem]{Lemma}
\newtheorem{defin}{Definition}[section]

\newtheorem{exa}{Example}[section]
\newenvironment{example}{\begin{exa}\rm}{\end{exa}}
\newtheorem{theorem}[defin]{Theorem}

\newtheorem{proposition}[defin]{Proposition}

\newtheorem{qu}{Question}

\newenvironment{proof}{\noindent{\it Proof.}}{\hfill $\Box$\par\vspace{2.5mm}}
\newenvironment{remark}{\par\vspace{2.5mm}\noindent{\it Remark.}}{\par\vspace{2.5mm}}

\numberwithin{equation}{section}

\title{{\sc Sharp forms of Nevanlinna error terms in differential
equations}}
\author{{\sc Risto Korhonen}}
\date{}

\begin{document}

\renewcommand{\thelettertheorem}{\Alph{lettertheorem}}

\maketitle

\begin{abstract}
Sharp versions of some classical results in differential equations
are given. Main results consists of a Clunie and a Mohon'ko type
theorems, both with sharp forms of error terms. The sharpness of
these results is discussed and some applications to nonlinear
differential equations are given in the conluding remarks.

Moreover, a short introduction on the connection between
Nevanlinna theory and number theory, as well as on their relation
to differential equations, is given. In addition, a brief review
on the recent developments in the field of sharp error term
analysis is presented.
\end{abstract}

\thispagestyle{empty}

\maketitle

\section{Introduction}

\renewcommand{\thefootnote}{}
\footnotetext[1]{\emph{Mathematics Subject Classification 2000:}
30D35, 34M05, 11J97.}

Rolf Nevanlinna's theory of value distribution is undoubtedly one
of the great mathematical discoveries of the twentieth century.
Nevanlinna laid the foundations of the theory in a 100 pages
article appeared in 1925 \cite{nevanlinna}. This remarkable
contribution was later on described by Hermann Weyl as
\begin{quote}
\textsl{...one of the few great mathematical events in our century} \cite{weyl}.
\end{quote}
The only remaining big open question, proposed and partially
solved by Nevanlinna himself, was for a long time the
\textsl{Nevanlinna inverse problem}, which is essentially a
problem of finding a meromorphic function with prescribed
deficient values. This problem was solved in 1976 by David Drasin,
who used quasiconformal mappings to construct the desired
function~\cite{drasin}. Undoubtedly this result was a substantial
contribution to the theory, and some would even go so far as to
say that Drasin's work finally completed Nevanlinna theory. But,
no~doubt, the richness of value distribution theory goes much
further than that. Not only it has numerous applications in the
fields of differential and functional equations, but there is a
profound relation between Nevanlinna theory and number theory,
which also extends to the theory of differential equations. In
this paper we will concentrate on this deep connection between
these three theories.

The present paper is organized as follows. We start by recalling
the necessary notation in Section \ref{prelim}. We then continue
by giving a short introduction to the deep connection between
Nevanlinna theory and number theory in Section \ref{nevnum}. This
is by no means complete review of the topic, and, therefore we
refer to \cite{aimo1}, \cite{yech}, \cite{ru} and \cite{vojta} for
a more comprehensive survey on the connection between there two
theories. Then, in Section \ref{errort}, we give a brief review on
the recent developments in the sharp error term analysis and
recall how this research is connected to the relation between
Nevanlinna theory and number theory. In Section~\ref{diffsection}
we discuss how the theory of differential equations is related to
the above topics. We also recall an important auxiliary result
which we will use to prove our main results in
Section~\ref{proofs}. Sections~\ref{results}--\ref{discuss}  are
the main sections of the present paper. They include the statement
of our results in Section~\ref{results}, the proofs of theorems in
Section~\ref{proofs} and concluding remarks in
Section~\ref{discuss}.

\section{Prerequisites}\label{prelim}

We use the standard notation of Nevanlinna theory, in other words
$m(r,f)$, $N(r,f)$ and $T(r,f)$ denote the \textsl{proximity
function}, the \textsl{counting function} and the
\textsl{characteristic function} of $f$, respectively. See for
example \cite{Hayman} and \cite{Laine} for explicit definitions
and basic results of Nevanlinna theory.

Consider an algebraic differential equation
    \begin{equation} \label{algdiii}
    P\left(z,w,w',\ldots,w^{(n)}\right) = 0
    \end{equation}
with meromorphic coefficients. We may write (\ref{algdiii}) in the form
    \begin{equation} \label{algdifeq2}
    \sum_{\lambda=(j_0,\ldots,j_n)\in I}
    a_{\lambda}(z)w^{j_0}(w')^{j_1}\cdots(w^{(n)})^{j_n} = 0,
    \end{equation}
where $I$ is a finite set of \textsl{multi-indices} $\lambda$, and
$a_{\lambda}$  are meromorphic functions for all $\lambda\in I$.
While considering equation \eqref{algdifeq2}, coefficients
$a_{\lambda}$ are often chosen to satisfy certain growth
conditions, for instance $T(r,a_{\lambda})=S(r,f)$ for $\lambda\in
I$. This will not, however, be necessary for the main results of
the present paper. The \textsl{degree} of a single term in
(\ref{algdifeq2}) is defined by
    \begin{equation*}\label{aste}
    d(\lambda):=j_0 + j_1 + \cdots + j_n,
    \end{equation*}
and its \textsl{weight} by
    \begin{equation*}
    w(\lambda):=j_1 + 2j_2 + \cdots + nj_n.
    \end{equation*}
Furthermore, we define the degree, $d(P)$, and the weight, $w(P)$, of the differential polynomial
$P$ as the maximal degree and the maximal weight of all terms of $P$,
respectively. Moreover, given a multi-index $\lambda=(j_0,\ldots,j_n)$, we denote
    \begin{equation*}
    |\lambda|:=j_0+\ldots+j_n.
    \end{equation*}
Finally, if a set $A$ satisfies
$A\subset \sqbrac{0,2\pi}$, then $\abs{A} := \int_A
\frac{dt}{2\pi}$.

The following, slightly modified version of the classical Jensen
inequality is needed to prove an important auxiliary result below.

\begin{lettertheorem}[Jensen]\label{jensen}
Let $(X,\mu)$ be a measure space such that $\mu(X)=\delta$, where
$0<\delta\leq 1$. Let $\chi$ be a concave function on the interval
$(a,+\infty)$, where $a\geq -\infty$. Then, for any integrable
real valued function $f$ with $a<f(x)<+\infty$ for almost all $x$
in $X$, we have
    \begin{equation*}
    \int_X \chi(f(x))d\mu(x) \leq \delta\cdot\chi\brac{\frac{1}{\delta}\int_X f(x) d\mu(x)}.
    \end{equation*}
\end{lettertheorem}

For a detailed proof of this refinement, see \cite[Theorem 4.1]{bk}.

Theorem \ref{jensen} has the following useful corollary, which we
name here as lemma. The proof is almost identical to that of
\cite[Corollary 4.2]{bk}, but since the assertion is somewhat more
general, we include the details of the proof here for convenience.

\begin{letterlemma}\label{mlemma}
Let $f(z)$ be a meromorphic function in $\set{z:\abs{z}<R}$, where
$0<R\leq\infty$, let $H\subset [0,2\pi]$ and let $0<\alpha<1$. Then for $0<r<R$,
    \begin{equation*}
     \int_H \log^{+} |f(re^{i\theta})|\,\frac{d\theta}{2\pi} \leq \frac{1}{\alpha} \brac{ \log^{+} \int_H
    \abs{f(re^{i\theta})}^{\alpha}\,
    \frac{d\theta}{2\pi} + e^{-1}}.
    \end{equation*}
\end{letterlemma}

\begin{proof}
First define
    \begin{equation*}
    E := \set{\theta\in H : \abs{f(re^{i\theta})} > 1},
    \end{equation*}
and $F := H \backslash E$. Then, by Theorem
\ref{jensen},
    \begin{eqnarray*}
     \int_H \log^{+} |f(re^{i\theta})|\,\frac{d\theta}{2\pi} &=& \frac{1}{\alpha}\int_H \log^{+}
    \abs{f(re^{i\theta})}^{\alpha}\frac{d\theta}{2\pi}
    \\ &=& \frac{1}{\alpha}\int_E \log
    \abs{f(re^{i\theta})}^{\alpha}\frac{d\theta}{2\pi}
    + \frac{1}{\alpha}\int_F \log^{+}
    \abs{f(re^{i\theta})}^{\alpha}\frac{d\theta}{2\pi}
    \\ & \leq & \frac{\abs{E}}{2\pi\alpha}\log\brac{\frac{2\pi}{\abs{E}}\int_E
    \abs{f(re^{i\theta})}^{\alpha}\frac{d\theta}{2\pi}}\\
    & \leq & \frac{1}{\alpha}\brac{\log^{+}\brac{\int_H
    \abs{f(re^{i\theta})}^{\alpha}\frac{d\theta}{2\pi}} +
    \frac{\abs{E}}{2\pi}\log\frac{2\pi}{\abs{E}}}.
    \end{eqnarray*}
Since $x\log(x^{-1})\leq e^{-1}$ for all $0<x<1$, the assertion
follows.
\end{proof}

\section{A connection between Nevanlinna theory and number theory} \label{nevnum}

What could be described as the core of Nevanlinna theory is the
Second main theorem, which is a deep result concerning the value
distribution of meromorphic functions. It contains Picard's
theorem,
\begin{quote}
\textsl{a transcendental meromorphic function assumes infinitely often all values in the plane except at most two} \cite{rudin},
\end{quote}
as a special case, and, in its general form, is stated as follows
\cite{nevanlinna}:
\begin{lettertheorem}[Nevanlinna]
Let $f$ be a non-constant meromorphic function, let $q\geq 2$ and let $a_1,\ldots,a_q\in\C$ be distinct points. Then
    \begin{equation*}
    m(r,f)+\sum_{j=1}^q m\left(r,\frac{1}{f-a_j}\right) \leq
        2T(r,f)-N_1(r,f)+S(r,f),
    \end{equation*}
where $N_1(r,f)$ is positive and given by
    \begin{equation*}
    N_1(r,f)=N\left(r,\frac{1}{f'}\right)+2N(r,f)-N(r,f').
    \end{equation*}
 \end{lettertheorem}

We now recall another classical and deep result, this time from
number theory. The proof of the following Thue-Siegel-Roth theorem
\cite{roth} yielded Klaus Roth a Fields medal.
\begin{lettertheorem}[Roth]
Let $\alpha$ be an algebraic irrational number, and let
$\varepsilon>0$. Then there are only finitely many solutions $p/q$
to
    \begin{equation*}
    \left|\alpha -\frac{p}{q}   \right| \leq \frac{1}{q^{2+\varepsilon}},
    \end{equation*}
where $p$ and $q$ are relatively prime integers.
\end{lettertheorem}

\begin{table}[h]
\begin{center}
\caption{Vojta's Dictionary}\label{vojtatab}
{\small
\begin{tabular}{p{7.8cm}|p{7.8cm}}
\textbf{Value Distribution} & \textbf{Diophantine Approximation}\\
\hline
A non-constant meromorphic function & An infinite $\{x\}$ in a number field $F$ \\
A radius $r$ & An element $x$ of $F$\\
A finite measure set $E$ of radii & A finite subset of $\{x\}$\\
An angle $\theta$ & An embedding $\sigma:F\hookrightarrow \C$\\
The modulus $|f(re^{i\theta})|$ & The modulus $|x|_{\sigma}=|\sigma(x)|$\\
The order of a pole: $\ord_z f$ & The order of $x$ in the prime ideal: $\ord_{\wp}x$\\
$\log\frac{r}{|z|}$ & $[F_{\wp}:\Q_p]\log p$\\
Characteristic function:  \newline $\displaystyle T(r,f)=m(r,f)+N(r,f)$ & Logarithmic height: \newline $\displaystyle h(x)=\frac{1}{[F:\Q]}\sum_{\sigma:F\hookrightarrow \C} \log^{+}|x|_{\sigma}+N(x,F)$\\
Mean proximity function:  \newline $\displaystyle m\left(r,\frac{1}{f-a}\right)=\int_0^{2\pi}\log^{+}\left|\frac{1}{f(re^{i\theta})-a}\right|\frac{d\theta}{2\pi}$ & Mean proximity function:  \newline $\displaystyle m\left(x,\frac{1}{F-a}\right)=\sum_{\sigma:F\hookrightarrow \C} \log^{+}\left|\frac{1}{x-a}\right|_{\sigma}$ \\
Counting function:  \newline $\displaystyle N\left(r,\frac{1}{f-a}\right)$\newline $\displaystyle{}\quad =\sum_{0<|z|<r}\ord_z^{+}(f-a)\log\frac{r}{|z|}+\ord_0^{+}\log r $ & Counting function:  \newline $\displaystyle N\left(x,\frac{1}{F-a}\right)$ \newline $\displaystyle {}\quad= \frac{1}{[F:\Q]}\sum_{\wp\in F} \ord_{\wp}^{+}(x-a)[F_{\wp}:\Q_p]\log p$ \\
First main theorem:  \newline $\displaystyle T(r,f)=T\left(r,\frac{1}{f-a}\right) + O(1)$ & Height property:  \newline $\displaystyle N\left(x,\frac{1}{F-a}\right) + m\left(x,\frac{1}{F-a}\right) = h(x) + O(1)$ \\
Weaker Second main theorem:  \newline $\displaystyle (q-2)T(r,f)-\sum_{j=1}^q N\left(r,\frac{1}{f-a_j}\right) \leq \veps T(r,f)$ & Roth's theorem:  \newline $\displaystyle (q-2)h(x) -\sum_{j=1}^q N\left(x,\frac{1}{F-a_j}\right) \leq \veps h(x) $ \\
\end{tabular}
}
\end{center}
\end{table}

At the first glance Second main theorem and Roth's theorem, like
Nevanlinna theory and number theory in general, seem to have very
little, or nothing, to do with each other. However, Charles Osgood
observed that this is quite not the case. He realized that the
occurrence of the number two in both of these classical results is
by no means an accident. Osgood was the first to draw attention to
the formal connections between Nevanlinna theory and number
theory, Roth's theorem and Second main theorem being just one
example~\cite{osgood}. Paul Vojta made this connection deeper and
more extensive by introducing so called \textsl{Vojta's
dictionary}, where he associated the terminology in Nevanlinna
theory and Diophantine approximation theory in an explicit manner
\cite{vojta}. This approach is applied also in Ru's monograph
\cite{ru}, which offers extensive correspondence tables, in
addition to Vojta's dictionary, between results and conjectures in
Nevanlinna theory and Diophantine approximation theory.  We are
settled with giving a few passages from Vojta's dictionary in
Table \ref{vojtatab}, just to briefly illustrate the basic idea.
For the complete representation, including a description of the
relevant terminology, we refer to \cite{vojta} and \cite{ru}.

What still remains as an `ultimate challenge', which would really
complete the analogue between these two fields of mathematics, is
to prove results in number theory, for instance Roth's theorem, by
using methods from Nevanlinna theory. The aim is roughly the
following: First take a result in Nevanlinna theory, then switch
the corresponding terms of Nevanlinna theory with the concepts in
number theory by using Vojta's dictionary, and finally follow the
reasoning in the proof of the original result. The result of this
procedure would be a theorem in number theory. This description
is, of course, a crude simplification of the real situation.
However, since Nevanlinna theory seems to be somewhat ahead of
number theory in the sense of the correspondence described above,
this approach would yield strong results in number theory, which
so far remain unproved, including the famous $abc$-conjecture due
to Masser and Oesterl\'e. Unfortunately, the dictionary does not
extend to all concepts needed to follow the above algorithm. In
particular, no counterpart for the derivative of a meromorphic
function has yet been discovered. Therefore it is, so far,
impossible to translate any proofs from Nevanlinna theory to
number theory.

\section{Sharp error term analysis}\label{errort}

Serge Lang observed in the beginning of 1990's that the analogy
between Nevanlinna theory and Diophantine approximation theory
motivates the search for a sharp form of the error term in
Nevanlinna's main theorems \cite{lang}. Lang's perception led to
extensive research activity.  Joseph Miles  obtained a sharp form
of the lemma on the logarithmic derivative both in the complex
plane and in the unit disc \cite{miles}. He also showed that his
results were essentially the best possible. His results were
studied further by Marcus Jankowski \cite{jankowski}, who showed
that a slightly smaller error term may be obtained, but with a
cost of a larger exceptional set.  Moreover, Aimo Hinkkanen proved
\cite{hinkkanen} a version of the Second main theorem where the
error term was sharp in the sense of \cite{miles}.

One of the fundamental results used in the sharp error term
analysis is the following estimate by Gol'dberg and Grinshtein \cite{GG}:
        \begin{equation}\label{ggestim}
        m\brac{r,\frac{f'}{f}} \leq \log^{+} \brac{\frac{T(\rho,f)}{r}
        \frac{\rho}{\rho-r}} + 5.8501,
        \end{equation}
where $0<r<\rho<R$ and $f$ is a meromorphic function in the disc
$|z|<R$ such that $f(0)=1$. Inequality \eqref{ggestim} has been
used, in one way or another, in most papers involving sharp error
term analysis. The constant term $5.8501$ was later on improved to
$5.3078$ by Djamel Benbourenane and the author   \cite{bk}.
Moreover, inequality \eqref{ggestim} was generalized to higher
order derivatives in a paper by Janne Heittokangas et
al.~\cite{hkr}, where one of the main results was the following:

\begin{lettertheorem}\label{mainthm}
Let $f$ be a meromorphic function in $\C$ such that $f^{(j)}$ does
not vanish identically, and let  $k$ and $j$ be integers such that
$k>j\geq 0$. Then there exists an $r_0>1$ such that
    \begin{equation}\label{mainres}
        m\brac{r,\frac{f^{(k)}}{f^{(j)}}} \leq (k-j)\log^{+}\brac{
        \frac{T(\rho,f)}{r}
        \frac{\rho}{\rho-r}} + \log\frac{k!}{j!} + (k-j)5.3078
    \end{equation}
for all $r_0<r<\rho<\infty$.
\end{lettertheorem}

Analogous result was shown to be true also in the unit disc. It
was also shown in \cite{hkr} that the best possible constant in
inequality \eqref{mainres} is greater or equal to
$(k-j)\log\frac{\pi}{e}$. This obviously gives also a limit to how
much the constant in \eqref{ggestim} can be further improved.
These sharpness considerations were accomplished by means of a
surprisingly simple auxiliary function
$$
f(z)=e^{z^n},
$$
which satisfies
 \begin{equation*}
         m\brac{r,\frac{f^{(k)}}{f^{(j)}}} \geq (k-j)\log^{+}\brac{
        \frac{T(\rho,f)}{r}\frac{\rho}{\rho-r}} +
        (k-j)\log\frac{\pi}{e}-\veps(r,\rho)
    \end{equation*}
for all sufficiently large $r$, where $\rho$ is set to be
$\rho=\frac{n}{n-1}r$, $n\in\N\setminus\{1\}$. The expression
$\veps(r,\rho)$ tends to zero as  $r$ and $n$ both approach to
infinity. We summarize the above observations as follows:

\begin{proposition}
Let $0<r<\rho<R$. The smallest possible constant $\kappa$, such that
        \begin{equation}\label{ggestim2}
        m\brac{r,\frac{f'}{f}} \leq \log^{+} \brac{\frac{T(\rho,f)}{r}
        \frac{\rho}{\rho-r}} + \kappa
        \end{equation}
for all meromorphic functions $f$ in the disc $|z|<R$, satisfies
$$0.1447 <  \log\frac{\pi}{e} \leq \kappa < 5.3078.$$
\end{proposition}

It is somewhat surprising that the small constant $\kappa$ cannot
be absorbed into the other, much larger term on the right hand
side of inequality \eqref{ggestim2}. It would be interesting to
see what the reason for this is. It is also possible that this
fact has a counterpart in number theory as well. But first, one
would need to find what the best possible constant is. This,
however, may prove to be a difficult task, and, since the
(partially numerical) methods used in \cite{GG} and \cite{bk} are
quite technical and complicated, it is likely that completely new
ideas are needed to find the desired constant.

\section{Differential equations}\label{diffsection}

Another connection between Nevanlinna theory and number theory can
be found through the theory of differential equations. On one
hand, Nevanlinna theory offers an efficient way to study
meromorphic solutions of differential equations, since it allows
the study of the properties of a solution without knowing its
explicit form. A good example of this efficacy is Kosaku Yosida's
simple alternative proof for classical Malmquist's theorem
\cite{yosida}. On the other hand, differential equations may be
used as a tool to study transcendental numbers. This approach has
been used, for instance, by Andrei Shidlovskii in his
monograph~\cite{shidlovskii}. Therefore differential equations can
be understood as a `bridge' between Nevanlinna theory and
transcendental number theory.

We now state two further examples from differential equations,
where Nevanlinna theory has been successfully utilized to obtain
accurate information about the meromorphic solutions. We start
with the  Clunie lemma \cite{clunie}, which is an efficient tool
to analyze the density of poles of solutions of certain kind of
algebraic differential equations.

\begin{letterlemma}[Clunie]\label{finorclunie}
Let $f$ be a transcendental meromorphic solution
of the differential equation
    \begin{equation*}
    f^n P(z,f) = Q(z,f),
    \end{equation*}
where $P(z,f)$ and $Q(z,f)$ are polynomials in $f$ and its
derivatives with meromorphic coefficients, say
$\{a_{\lambda}:\lambda\in I\}$, such that
$m(r,a_{\lambda})=S(r,f)$ for all $\lambda\in I$. If the total
degree of $Q(z,f)$ as a polynomial in $f$ and its derivatives is
$\leq n$, then
    \begin{equation*}
    m\brac{r,P(z,f)} = S(r,f).
    \end{equation*}
\end{letterlemma}

The second result is the following Mohon'ko's theorem
\cite[Theorem 6]{mohonko}. It can be used, similarly as Clunie's
lemma may be used to study pole distribution, to study value
distribution of meromorphic solutions of differential equations.

\begin{lettertheorem}[Mohon'ko]\label{finormok}
Let $P(z,u)$ be a polynomial in $u$ and its derivatives with
meromorphic coefficients, say $\{a_{\lambda}:\lambda\in I\}$, such
that $m(r,a_{\lambda})=S(r,f)$ for all $\lambda\in I$, and let
$u=f$ be a transcendental meromorphic solution of the differential
equation
$$
P(z,u)=0.
$$
If $P(z,0) \not\equiv 0$, then
$$
m\left(r,\frac{1}{f}\right)=S(r,f).
$$
\end{lettertheorem}

The fact that connects these two classical results is that both of
their proofs are based on a natural generalization
    $$
    m\left(r,\frac{f^{(n)}}{f}\right) = S(r,f)
    $$
of the lemma on the logarithmic derivative. Therefore, one might
expect that inequality~\eqref{mainres} would yield even stronger
results. Indeed this is true, but to obtain more natural results
somewhat different auxiliary lemma is needed, although it is in
the same spirit as inequality~\eqref{mainres}. The following lemma
was proved in \cite{hkr}.

\begin{letterlemma}\label{intlemma}
Let $f$ be a non-polynomial meromorphic function in $\C$. Let $k$ and $j$ be integers satisfying $k>j\geq 0$,
$\alpha$ and $\beta$ be constants satisfying $0<\alpha(k-j)<1$
and $0<\beta< 1$, and let $\varepsilon >0$. Then there exists an
$r_0(\varepsilon)>1$ such that, for all $r_0<r<\rho<\infty$,
        \begin{equation*}
        \int_0^{2\pi}
        \abs{\frac{f^{(k)}(re^{i\theta})}{f^{(j)}(re^{i\theta})}}^{\alpha}
        \frac{d\theta}{2\pi} \leq
        \left(C\big(\alpha(k-j),\beta\big)+\varepsilon\right)\brac{\frac{k!}{j!}}^{\alpha}
        \brac{\frac{T(\rho,f)}{r}\frac{\rho}{\rho-r}}^{\alpha(k-j)},
        \end{equation*}
where
        \begin{equation}\label{constant-c}
        \begin{split}
        C(a,b) = & \brac{\frac{2}{1-b}}^a +
        \frac{\sec\brac{\frac{a\pi}{2}}}{b^a} \brac{4+\left(
        2^{\frac{1+a}{1-a}}+2^{\frac{2+a}{1-a}}
        \right)^{1-a}}. \\
        \end{split}
        \end{equation}
\end{letterlemma}

\begin{remark}
 It was also calculated in
 \cite{bk} that
    \begin{equation*}
        \frac{\log (C(\alpha,\beta) + \varepsilon) + e^{-1}}{\alpha} +\veps< 5.3078,
    \end{equation*}
when $\alpha=0.815508$,
$\beta=0.845890$, and   $\varepsilon >0$ is chosen to be sufficiently small.
\end{remark}

Lemma \ref{intlemma} and the remark above will be repeatedly applied in the proofs of our main results.

\section{Statement of results}\label{results}

Understanding the whole structure behind the connections between
these three seemingly disconnected fields of mathematics,
differential equations, Nevanlinna theory and number theory, could
be described as assembling a complicated puzzle without the help
of an exemplar. Some of the pieces are already in place, and so we
are beginning to see how the big picture looks like. However,
large part of the work is still undone, especially in
understanding how differential equations fit into this picture.

We will make here, best to our knowledge, the first attempt to
gain more information about these connections by close examination
of error terms in differential equations. We give versions of the
two classical results stated in Section \ref{diffsection} with
sharp forms of error terms. We start with a Clunie type result.
Although Clunie modestly called his original result lemma, we
label ours as a theorem since it is one of the main results of the
present paper.

\begin{theorem}\label{clunie}
Let $f$ be a transcendental meromorphic solution
of the differential equation
    \begin{equation*}
    f^n P(z,f) = Q(z,f),
    \end{equation*}
where
    \begin{equation*}
    P(z,f)= \sum_{\lambda\in I} a_{\lambda} f^{i_0}(f')^{i_1}\cdots(f^{(m)})^{i_m}
    \end{equation*}
and
    \begin{equation*}
    Q(z,f)= \sum_{\mu\in J} b_{\mu} f^{j_0}(f')^{j_1}\cdots(f^{(k)})^{j_k}.
    \end{equation*}
If the total degree of $Q(z,f)$ as a polynomial in $f$ and its derivatives is $\leq n$,
then there exists $r_0$ such that
    \begin{equation*}\label{cin}
    \begin{split}
    m\brac{r,P(z,f)} &\leq \big(w(P)+w(Q)\big) \log^{+} \brac{\frac{T(\rho,f)}{r}\frac{\rho}{\rho-r}}
    +  \max_{\lambda\in I}
    m(r,a_{\lambda}) \\ &\quad + \max_{\mu\in J}     m(r,b_{\mu})
    + d(P)\log m! + d(Q)\log k! + \log^{+} \card(I)  \\ &\quad+ \log^{+} \card(J)
       +  \big(w(P)+w(Q)\big)5.3078\\
    \end{split}
    \end{equation*}
for all $r_0<r<\rho<\infty$.
\end{theorem}

The next result is a sharp form of Mohon'ko's theorem.

\begin{theorem}\label{mokh}
Let $f$ be a transcendental meromorphic solution
of the differential equation
    \begin{equation}\label{peqn}
    P(z,f) = 0,
    \end{equation}
where
    \begin{equation*}
    P(z,f)= \sum_{\lambda\in I} a_{\lambda} f^{i_0}(f')^{i_1}\cdots(f^{(m)})^{i_m}.
    \end{equation*}
If $a_0:=P(z,0)\not\equiv 0$, then there exists $r_0$ such that
    \begin{equation*}\label{min}
    \begin{split}
    m\brac{r,\frac{1}{f}} &\leq w(P) \log^{+} \brac{\frac{T(\rho,f)}{r}\frac{\rho}{\rho-r}} +
    \max_{\lambda\in I \atop |\lambda|\geq 1}
    m(r,a_{\lambda}) + m\left(r,\frac{1}{a_0}\right)  \\ &\quad + d(P)\log m!
        + \log^{+} (\card(I)-1)+  w(P)5.3078\\
    \end{split}
    \end{equation*}
for all $r_0<r<\rho<\infty$.
\end{theorem}

\section{Proofs of theorems}\label{proofs}

The basic ideas behind the proofs of our results are adopted from
the original works of Clunie and Mohon'ko. These ideas are
complemented by more recent results taking the form of the error
term into account.

\subsection*{The proof of Theorem \ref{clunie}}

Defining
    \begin{equation}\label{sets}
    \begin{split}
    E_1 &=\left\{\varphi\in[0,2\pi) : |f(re^{i\varphi})|<1 \right\}\\
    E_2 &=[0,2\pi)\setminus E_1\\
    \end{split}
    \end{equation}
and denoting
    \begin{equation}\label{fixp}
    P(z):=P(z,f)=: \sum_{\lambda\in I} a_{\lambda} P_{\lambda}(z),
    \end{equation}
we consider the proximity function $m(r,P)$ in two parts:
    \begin{equation*}
    m(r,P) = \int_{E_1} \log^{+}|P(re^{i\varphi})|\,\frac{d\varphi}{2\pi} + \int_{E_2}
    \log^{+}|P(re^{i\varphi})|\,\frac{d\varphi}{2\pi}.
    \end{equation*}

We start by looking at $P_{\lambda}(z)$ in the set $E_1$. Let
$\alpha\in(0,\frac{1}{w(\lambda)})$ and
$\lambda=(i_0,i_1,\ldots,i_m)$. Then, by the H\"older inequality
and Lemma \ref{intlemma},  there exists $r_1$ such that
    \begin{equation}\label{monsteri1}
    \begin{split}
    \int_{E_1} |P_{\lambda}(re^{i\varphi})|^{\alpha}&\,\frac{d\varphi}{2\pi}  \leq  \int_0^{2\pi} \left|\frac{f'(re^{i\varphi})}{f(re^{i\varphi})}\right|^{\alpha i_1}
        \cdots
        \left|\frac{f^{(m)}(re^{i\varphi})}{f(re^{i\varphi})}\right|^{\alpha i_m}
        \frac{d\varphi}{2\pi}
        \\
        &\leq \left(\int_0^{2\pi}
        \left|\frac{f'(re^{i\varphi})}{f(re^{i\varphi})}\right|^{\alpha
        w(\lambda)}\frac{d\varphi}{2\pi}\right)^{\frac{i_1}{w(\lambda)}}
        \cdots
        \left(\int_0^{2\pi}
        \left|\frac{f^{(m)}(re^{i\varphi})}{f(re^{i\varphi})}\right|^{\frac{\alpha
        w(\lambda)}{m}}\frac{d\varphi}{2\pi}\right)^{\frac{m i_m}{w(\lambda)}}
        \\
        &\leq \left(    \left(C\big(\alpha w(\lambda),\beta\big)+\varepsilon\right)
         \brac{\frac{T(\rho,f)}{r}\frac{\rho}{\rho-r}}^{\alpha w(\lambda)}      \right)^{\frac{i_1}{w(\lambda)}}
         \\
        &\quad\cdots
        \left(  \left(C\big(\alpha w(\lambda),\beta\big)+\varepsilon\right)(m!)^{\frac{\alpha w(\lambda)}{m}}
        \brac{\frac{T(\rho,f)}{r}\frac{\rho}{\rho-r}}^{\alpha w(\lambda)} \right)^{\frac{m i_m}{w(\lambda)}}
        \\
        &= \left(C\big(\alpha w(\lambda),\beta\big)+\varepsilon\right)1\cdot (2!)^{\alpha i_2}\cdots (m!)^{\alpha i_m}
        \brac{\frac{T(\rho,f)}{r}\frac{\rho}{\rho-r}}^{\alpha w(\lambda)}
        \end{split}
    \end{equation}
for all $r_1<r<\rho<\infty$. Now, remembering notation \eqref{fixp}, we get
    \begin{equation}\label{vali}
    \begin{split}
    \int_{E_1} \log^{+}|P(re^{i\varphi})|\,\frac{d\varphi}{2\pi} &=
        \int_{E_1} \log^{+}\left|\sum_{\lambda\in I} a_{\lambda}(re^{i\varphi})
     P_{\lambda}(re^{i\varphi})  \right|\,\frac{d\varphi}{2\pi}
    \\
    &\leq \int_{E_1} \log^{+} \left(\card(I)\max_{\lambda\in I} \left|a_{\lambda}(re^{i\varphi})
     P_{\lambda}(re^{i\varphi})  \right|\right)\,\frac{d\varphi}{2\pi}
    \\
    &\leq  \max_{\lambda\in I}\left( \int_{E_1} \log^{+}\left| P_{\lambda}(re^{i\varphi})\right|\,\frac{d\varphi}{2\pi}
     + m(r,a_{\lambda})\right)
      +\log^{+} \card(I).
    \\
    \end{split}
    \end{equation}
Hence, by Lemma \ref{mlemma}, inequality  \eqref{monsteri1} and the remark below Lemma \ref{intlemma}, we obtain
    \begin{equation}\label{e1}
    \begin{split}
    \int_{E_1}& \log^{+}|P(re^{i\varphi})|  \,\frac{d\varphi}{2\pi}
    \leq \max_{\lambda\in I}\left(   \frac{1}{\alpha} \brac{ \log^{+} \int_{E_1}
       \abs{P_{\lambda}(re^{i\varphi})}^{\alpha}
        \frac{d\varphi}{2\pi} + e^{-1}} \right)
        \\
        &\quad+\max_{\lambda\in I} m(r,a_{\lambda})  +\log^{+} \card(I)
      \\
      &\leq \max_{\lambda\in I} \frac{1}{\alpha}  \log^{+}\brac{
       \left(C\big(\alpha w(\lambda),\beta\big)+\varepsilon\right)1\cdot (2!)^{\alpha i_2}\cdots (m!)^{\alpha i_m}
        \brac{\frac{T(\rho,f)}{r}\frac{\rho}{\rho-r}}^{\alpha w(\lambda)}}
        \\
       &\quad+ \frac{e^{-1}}{\alpha} + \max_{\lambda\in I}  m(r,a_{\lambda})
      +\log^{+} \card(I)
      \\
      &\leq w(P) \log^{+} \brac{\frac{T(\rho,f)}{r}\frac{\rho}{\rho-r}} +
       \max_{\lambda\in I}  m(r,a_{\lambda}) + d(P) \log^{+} m!  +\log^{+} \card(I)
        \\
       &\quad+  w(P) \max_{\lambda\in I}  \frac{\log^{+}\left(C\big(\alpha w(\lambda),\beta\big)+\varepsilon\right)
       + e^{-1}}{\alpha w(\lambda)}
          \\
      &\leq w(P) \log^{+} \brac{\frac{T(\rho,f)}{r}\frac{\rho}{\rho-r}} +  \max_{\lambda\in I}  m(r,a_{\lambda}) + d(P)
      \log^{+} m!  +\log^{+} \card(I)
      \\
       &\quad +  w(P)5.3078
      \\
    \end{split}
    \end{equation}
for all $r_1<r<\rho<\infty$.

To consider the set $E_2$, recall that we assumed $j_0 + \cdots +
j_k\leq n $  for all $\mu=(j_0, \cdots,  j_k)\in J$. Therefore,
    \begin{equation}\label{e2}
    \begin{split}
    \int_{E_2} & \log^{+}|P(re^{i\varphi})|\,\frac{d\varphi}{2\pi} =
        \int_{E_2} \log^{+}\left|\frac{1}{f(re^{i\varphi})^n}\sum_{\mu\in J}
    b_{\mu} (re^{i\varphi})
      f'(re^{i\varphi})^{j_1}
        \cdots
        f^{(m)}(re^{i\varphi})^{j_k}
       \right|\,\frac{d\varphi}{2\pi}
    \\
    &\leq   \int_{E_2} \log^{+}\left|\sum_{\mu\in J} b_{\mu} (re^{i\varphi})
      \left(\frac{f'(re^{i\varphi})}{f(re^{i\varphi})}\right)^{j_1}
        \cdots
              \left(\frac{f^{(k)}(re^{i\varphi})}{f(re^{i\varphi})}\right)^{j_k}
       \right|\,\frac{d\varphi}{2\pi}
    \\
    &\leq  \max_{\mu \in J}\left( \int_{E_2} \log^{+}\left|
       \left(\frac{f'(re^{i\varphi})}{f(re^{i\varphi})}\right)^{j_1}
        \cdots
              \left(\frac{f^{(k)}(re^{i\varphi})}{f(re^{i\varphi})}\right)^{j_k}
                 \right|\,\frac{d\varphi}{2\pi}
     + m(r,b_{\mu})\right) \\ &\quad
      +\log^{+} \card(J).
    \\
    \end{split}
    \end{equation}
By continuing with essentially identical reasoning as in
inequalities \eqref{monsteri1}, \eqref{vali} and \eqref{e1}, we
finally conclude that there is a constant $r_2$ such that
    \begin{equation}\label{e3}
    \begin{split}
    \int_{E_2}  \log^{+}|P(re^{i\varphi})|\,\frac{d\varphi}{2\pi} & \leq w(Q) \log^{+}
    \brac{\frac{T(\rho,f)}{r}\frac{\rho}{\rho-r}} +  \max_{\mu\in J}  m(r,b_{\mu}) \\ &\quad + d(Q)
    \log^{+} k!  +\log^{+} \card(J)  +  w(Q)5.3078
      \end{split}
      \end{equation}
for all $r_2<r<\rho<\infty$.

The assertion follows by setting $r_0=\max\{r_1,r_2\}$, and by
combining estimates \eqref{e1} and~\eqref{e3}. \hfill $\Box$

\subsection*{Proof of Theorem \ref{mokh}}

We may write equation \eqref{peqn} in the form
    \begin{equation*}
    P(z,f)=a_{0}(z)+Q(z)=0,
    \end{equation*}
where $a_0(z)=P(z,0)\not\equiv 0$ by assumption, and where
    \begin{equation*}
    Q(z)= \sum_{\lambda\in I \atop |\lambda|\geq 1} a_{\lambda}
    f^{i_0}(f')^{i_1}\cdots(f^{(m)})^{i_m}.
    \end{equation*}
Recalling \eqref{sets}, we consider the situation again in the two
sets $E_1$ and $E_2$ separately. Clearly the integral $m(r,1/f)$
vanishes on  $E_2$.  On the other hand, in $E_1$,
    \begin{equation*}
    \frac{1}{|f|} \left| f^{i_0}(f')^{i_1}\cdots(f^{(m)})^{i_m}  \right| \leq \left|\frac{f'}{f}\right|^{i_1} \cdots
    \left|\frac{f^{(m)}}{f}\right|^{i_m}
    \end{equation*}
for each term of $Q(z)$. Therefore,
    \begin{equation*}
    \begin{split}
    m\left(r,\frac{1}{f}\right) &= \int_{E_1}
      \log^{+}\left|\frac{1}{f(re^{i\varphi})}\right|\,\frac{d\varphi}{2\pi}
    \\
    &= \int_{E_1} \log^{+}\left|\frac{a_0(re^{i\varphi})}{f(re^{i\varphi})}
    \frac{1}{a_0(re^{i\varphi})}\right|\,\frac{d\varphi}{2\pi}
    \\
    &= \int_{E_1} \log^{+}\left|\frac{Q(re^{i\varphi})}{f(re^{i\varphi})}
    \frac{1}{a_0(re^{i\varphi})}\right|\,\frac{d\varphi}{2\pi}
    \\
    &\leq \int_{E_1} \log^{+}\left|\frac{Q(re^{i\varphi})}{f(re^{i\varphi})}
    \right|\,\frac{d\varphi}{2\pi}  + m\left(r,\frac{1}{a_0}\right)
    \\
    &\leq   \int_{E_1} \log^{+}\left((\card(I)-1) \max_{\lambda\in I \atop |\lambda|\geq 1} \left|
    \frac{a_{\lambda}(re^{i\varphi})
    f(re^{i\varphi})^{i_0}f'(re^{i\varphi})^{i_1}\cdots f^{(m)}(re^{i\varphi})^{i_m}}{f(re^{i\varphi})}
    \right|\right)\,\frac{d\varphi}{2\pi}
    \\
     &\quad + m\left(r,\frac{1}{a_0}\right)
    \\
    &\leq  \max_{\lambda\in I \atop |\lambda|\geq 1}  \int_{E_1} \log^{+} \left|  \frac{
    f'(re^{i\varphi})}{f(re^{i\varphi})}\right|^{i_1}\cdots \left|  \frac{
    f^{(m)}(re^{i\varphi})}{f(re^{i\varphi})}\right|^{i_m}
    \,\frac{d\varphi}{2\pi} +\max_{\lambda\in I \atop |\lambda|\geq 1}  m(r,a_{\lambda})
    \\ &\quad+ m\left(r,\frac{1}{a_0}\right) + \log^{+}(\card(I)-1).
    \\
    \end{split}
    \end{equation*}
Again, by applying the reasoning used in inequalities
\eqref{monsteri1}, \eqref{vali} and \eqref{e1}, we conclude that
there is a constant $r_0$ such that
    \begin{equation*}
    \begin{split}
    m\left(r,\frac{1}{f}\right)&\leq w(P) \log^{+} \brac{\frac{T(\rho,f)}{r}\frac{\rho}{\rho-r}} +
    \max_{\lambda\in I \atop
    |\lambda|\geq 1}
    m(r,a_{\lambda})
    +m\left(r,\frac{1}{a_0}\right)  \\&\quad  + d(P) \log^{+} m!
    +\log^{+} (\card(I)-1) +  w(P)5.3078
      \\
      \end{split}
      \end{equation*}
for all $r_0<r<\rho<\infty$. \hfill $\Box$

\section{Concluding remarks}\label{discuss}

The obvious question which we have not addressed so far, is
whether or not our main estimates are, in one sense or another,
the best possible. We conclude the present paper by discussing
this question, and looking at some applications of our theorems.
Since Theorem \ref{clunie} can only be applied to nonlinear
differential equations, we start with the simplest nonlinear
equation, which is the Riccati equation.
\begin{example}\label{riccaex}
All meromorphic solutions of the Riccati differential equation
    \begin{equation}\label{ricca}
    w'=a(z)w^2+b(z)w+c(z),\tag{$R$}
    \end{equation}
with rational coefficients $a(z)$, $b(z)$ and $c(z)$, are of
finite order of growth, see \cite{Laine} for example. Assume that
$f$ is such a solution. Then, by Theorem \ref{clunie},
    \begin{equation*}
    \limsup_{r\longrightarrow +\infty} \frac{m(r,f)}{\log r} \leq \max\{0,\rho(f)-1\} + \di^{+}(a)
    + \max\{\di^{+}(b),\di^{+}(c)\},
    \end{equation*}
where, given a rational function $r=p/q$, where $p$ and $q$ are polynomials,
$$\di(r):=\deg(p)-\deg(q)$$
denotes the \textsl{degree} of a rational function $r$ \textsl{at infinity}, and
$$\di^{+}(r):=\max\{0,\di(r)\}.$$
Similarly, by Theorem \ref{mokh},
    \begin{equation*}
    \limsup_{r\longrightarrow +\infty} \frac{m\left(r,\frac{1}{f-q}\right)}{\log r} \leq \max\{0,\rho(f)-1\}
    + \max\left\{\di^{+}(a),\di^{+}(b)\right\}+ \di^{+}\left(\frac{1}{c}\right),
    \end{equation*}
where $q\in\C$. In particular, since $\tan z$ is a solution of
    \begin{equation*}
    w'=w^2+1,
    \end{equation*}
we have
    \begin{equation*}
    \limsup_{r\longrightarrow +\infty} \frac{m(r,\tan z)}{\log r}
    = \limsup_{r\longrightarrow +\infty} \frac{m\left(r,\frac{1}{\tan z}\right)}{\log r}=0,
    \end{equation*}
which is, of course, a well known fact.
\end{example}

In fact, Theorem \ref{clunie} correctly yields
    \begin{equation*}
    m(r,\tan z) = O(1),
    \end{equation*}
and similarly, by Theorem \ref{mokh}, we obtain
    \begin{equation*}
    m\left(r,\frac{1}{\tan z}\right) = O(1).
    \end{equation*}
Unfortunately, this is not enough to demonstrate the sharpness of
the main estimates. Indeed, since the weight of the Riccati
equation is one, it is not likely to yield any further
information. So, in order to get better examples, we must turn our
attention to higher order nonlinear differential equations, the
next obvious choice being the Painlev\'e differential equations.

\begin{example}\label{pl}
Consider the first, second and the fourth Painlev\'e differential equations
    \begin{equation}\label{I}
    w'' = 6w^2+z,\tag{$P_I$}
    \end{equation}
    \begin{equation}\label{II}
    w'' = 2w^3+zw+\alpha,\tag{$P_{II}$}
    \end{equation}
    \begin{equation}\label{IV}
    ww'' = \frac{1}{2}(w')^2 + \frac{3}{2}w^4 + 4zw^3 +
    2(z^2-\beta)w^2 + \gamma,\tag{$P_{IV}$}
    \end{equation}
where $\alpha,\beta,\gamma\in\C$. It is well known that all
solutions of these equations are meromorphic, and of finite order
of growth. It is also known that
    \begin{equation*}
    m(r,f)=O(\log r),
    \end{equation*}
if $f$ is a solution of any of the equations \eqref{I}, \eqref{II}
or \eqref{IV}. Now, by Theorem \ref{clunie}, we have in fact that
    \begin{equation}\label{es1}
    m(r,f_I)\leq 3\log r +O(1),
    \end{equation}
where $f_I$ is a solution of \eqref{I}. Similarly
    \begin{equation}\label{es2}
    m(r,f_{II})\leq 4\log r +O(1)
    \end{equation}
and
    \begin{equation}\label{es4}
    m(r,f_{IV})\leq 6\log r +O(1),
    \end{equation}
where $f_{II}$ and $f_{IV}$ are solutions of \eqref{II} and
\eqref{IV}, respectively. While deriving inequalities \eqref{es2},
\eqref{es2} and \eqref{es4}, we have used, in addition to Theorem
\ref{clunie}, the fact that $\rho(f_I)=5/2$, $\rho(f_{II})\leq 3$
and $\rho(f_{IV})\leq 4$, see \cite{shimomura} or
\cite{steinmetz}.
\end{example}

Finding out whether or not the estimates \eqref{es1}, \eqref{es2}
and \eqref{es4} are sharp, would possibly give an answer to the
question about sharpness of our main estimates. However, since
dealing with the Painlev\'e transcendents is highly complicated
issue in general, this is not an easy task, and, indeed is beyond
the scope of this paper. We are settled with leaving this matter
as an open question.

\bigskip

\noindent {\sc Department of Mathematical Sciences,
 Loughborough University, Loughborough, Leicestershire, LE11 3TU, UK.}

\smallskip

\noindent\emph{E-mail address:} {\tt r.j.korhonen@lboro.ac.uk}

\newpage\section*{Erratum:}

In the above article there is an error in the proofs of Theorem
6.1 and 6.2. The expression ``$\max_{\lambda\in I}$'' in
inequalities (7.4) and (7.5), and in the inequality on page 11,
should be replaced by ``$\sum_{\lambda\in I}$'', and similarly
``$\max_{\mu\in J}$'' in inequalities (7.6) and (7.7) should be
``$\sum_{\mu\in J}$''. This modification has a slight effect on
the final form of Theorems 6.1 and 6.2. The correct statements of
these theorems are as follows.

\bigskip

\noindent\textbf{Theorem 6.1.} Let $f$ be a transcendental
meromorphic solution of the differential equation
    \begin{equation*}
    f^n P(z,f) = Q(z,f),
    \end{equation*}
where
    \begin{equation*}
    P(z,f)=\sum_{\lambda\in I} P_\lambda(z,f)= \sum_{\lambda\in I} a_{\lambda} f^{i_0}(f')^{i_1}\cdots(f^{(m)})^{i_m}
    \end{equation*}
and
    \begin{equation*}
    Q(z,f)=\sum_{\mu\in J} Q_\lambda(z,f)= \sum_{\mu\in J} b_{\mu} f^{j_0}(f')^{j_1}\cdots(f^{(k)})^{j_k}.
    \end{equation*}
If the total degree of $Q(z,f)$ as a polynomial in $f$ and its
derivatives is $\leq n$, then there exists $r_0$ such that
    \begin{equation*}
    \begin{split}
    m\left(r,P(z,f)\right) &\leq \bigg(\sum_{\lambda\in I} w(P_\lambda)+\sum_{\mu\in J}w(Q_\mu)\bigg)
     \log^{+} \left(\frac{T(\rho,f)}{r}\frac{\rho}{\rho-r}\right) +  \sum_{\lambda\in I}
    m(r,a_{\lambda}) \\ &\quad + \sum_{\mu\in J}     m(r,b_{\mu})
    + \sum_{\lambda\in I} d(P_\lambda)\log m! + \sum_{\mu\in J} d(Q_\mu)\log k! + \log^{+} \mathrm{card}(I)
    \\ &\quad+ \log^{+} \mathrm{card}(J)
       +  \bigg(\sum_{\lambda\in I} w(P_\lambda)+\sum_{\mu\in J} w(Q_\mu)\bigg)5.3078\\
    \end{split}
    \end{equation*}
for all $r_0<r<\rho<\infty$.\bigskip

Here $d(P_\lambda)=i_0 + \cdots + i_m$ is the \textit{degree} of
$P_\lambda(z,f)=a_{\lambda}
f^{i_0}(f')^{i_1}\cdots(f^{(m)})^{i_m}$, and
$w(P_\lambda)=i_1+2i_2+ \cdots + mi_m$ denotes the \textit{weight}
of $P_\lambda(z,f)$.

\bigskip

\noindent\textbf{Theorem 6.2.} Let $f$ be a transcendental
meromorphic solution of the differential equation
    \begin{equation*}
    P(z,f) = 0,
    \end{equation*}
where
    \begin{equation*}
    P(z,f)= \sum_{\lambda\in I} P_\lambda(z,f) =\sum_{\lambda\in I} a_{\lambda} f^{i_0}(f')^{i_1}\cdots(f^{(m)})^{i_m}.
    \end{equation*}
If $a_0:=P(z,0)\not\equiv 0$, then there exists $r_0$ such that
    \begin{equation*}
    \begin{split}
    m\left(r,\frac{1}{f}\right) &\leq \sum_{\lambda\in I} w(P_\lambda)
    \log^{+} \left(\frac{T(\rho,f)}{r}\frac{\rho}{\rho-r}\right) +
    \sum_{\genfrac{}{}{0pt}{}{\lambda\in I}{|\lambda|\geq 1}}
    m(r,a_{\lambda}) + m\left(r,\frac{1}{a_0}\right)  \\ &\quad + \sum_{\lambda\in I} d(P_\lambda)\log m!
        + \log^{+} (\mathrm{card}(I)-1)+  \sum_{\lambda\in I} w(P_\lambda)5.3078\\
    \end{split}
    \end{equation*}
for all $r_0<r<\rho<\infty$.

\bigskip

Theorems 6.1 and 6.2 were applied to find upper bounds for the
proximity functions of the first, second and fourth Painlev\'e
transcendents. Due to the above mentioned error, these estimates
(8.1) -- (8.3) are also incorrect. The correct estimates are
    \begin{equation}
    m(r,f_I) \leq  4\log r +O(1) \tag{8.1},
    \end{equation}
    \begin{equation}
    m(r,f_{II}) \leq  5\log r +O(1) \tag{8.2},
    \end{equation}
    \begin{equation}
    m(r,f_{IV}) \leq  15\log r +O(1) \tag{8.3}.
    \end{equation}

\bigskip

The author wishes to express his thanks to those who pointed out
these errors to him.

\bigskip

\noindent Current address: {\sc University of Joensuu,
Mathematics, P.O. Box 111, FI-80101 Joensuu, Finland.}

\smallskip

\noindent\emph{E-mail address:} {\tt risto.korhonen@joensuu.fi}


\begin{thebibliography}{99}
\bibitem{aimo1} \textsc{A.~Hinkkanen}, `Sharp error term in the Nevanlinna theory',
Complex differential and functional equations: proceedings of the
summer school held in Mekrij\"arvi, July 30 -- August 3, 2000,
edited by Ilpo Laine (Report series, University of Joensuu,
Department of Mathematics, no.~5).
\bibitem{bk} \textsc{D.~Benbourenane} and \textsc{R.~Korhonen}, `On the growth of the
logarithmic derivative',   Comput.~Methods Funct.~Theory 1 (2002) 301--310.
\bibitem{yech} \textsc{W.~Cherry} and \textsc{Z.~Ye}, \textsl{Nevanlinna's theory of value
distribution. The second main theorem and its error terms} (Springer-Verlag, Berlin, 2001).
\bibitem{clunie} \textsc{J.~Clunie}, `On integral and meromorphic functions', J.~London Math.~Soc. 37 (1962) 17--27.
\bibitem{drasin} \textsc{D.~Drasin}, `The inverse problem of the Nevanlinna theory',  Acta Math.  138  (1976) 83--151.
\bibitem{GG} \textsc{A.~Gol'dberg} and \textsc{V.~Grinshtein}, `The logarithmic derivative
of a meromorphic function' (Russian), Mat. Zametki 19
(1976) 525--530 (English translation: Math. Notes 19
(1976) 320--323).
\bibitem{Hayman} \textsc{W.~K.~Hayman}, \textsl{Meromorphic Functions}
(Clarendon Press, Oxford, 1964).
\bibitem{hkr} \textsc{J.~Heittokangas}, \textsc{R.~Korhonen} and \textsc{J.~R\"atty\"a}, `Generalized logarithmic derivative estimates of Gol'dberg-Grinshtein type', to appear in Bull.~London Math.~Soc.
\bibitem{hinkkanen} \textsc{A.~Hinkkanen}, `A sharp form of Nevanlinna's
second fundamental theorem', Invent. Math. 108 (1992)
549--574.
\bibitem{jankowski} \textsc{M.~Jankowski}, `An estimate for the
logarithmic derivative of meromorphic functions', Analysis
14 (1994) 185--194.
\bibitem{Laine} \textsc{I.~Laine}, \textsl{Nevanlinna Theory and Complex
Differential Equations} (Walter de Gruyter, Berlin, 1993).
\bibitem{lang} \textsc{S.~Lang} and \textsc{W.~Cherry}, \textsl{Topics in Nevanlinna Theory},
Lecture Notes in Math.~1433 (Springer-Verlag, Berlin, 1990).
\bibitem{miles} \textsc{J.~Miles}, `A sharp form of the lemma on the logarithmic
derivative', J.~London Math.~Soc.~45 (1992) 243--254.
\bibitem{mohonko} \textsc{A.~Z. Mohon'ko} and \textsc{V.~D.~Mohon'ko}, `Estimates
for the Nevanlinna characteristic of some classes of meromorphic
functions and their applications to differential equations',
Sibirsk.~Mat.~Zh.~15 (1974) 1305--1322 [Russian]. Engl.~Trans.:
Siberian Math.~J.~15 (1974) 921--934.
\bibitem{nevanlinna} \textsc{R.~Nevanlinna}, `Zur Theorie der meromorphen Funktionen', Acta Math.~46 (1925) 1--99.
\bibitem{osgood} \textsc{C.~Osgood}, `Sometimes effective
Thue-Siegel-Roth-Schmidt-Nevanlinna bounds, or better', J.~Number
Theory 21 (1985) 347--398.
\bibitem{roth} \textsc{K.~F.~Roth}, `Rational approximations to algebraic numbers', Mathematika~2 (1955) 1--20.
\bibitem{ru} \textsc{M.~Ru}, \textsl{Nevanlinna theory and its relation to Diophantine approximation}
 (World Scientific, Singapore-New Jersey-London-Hong Kong, 2001).
\bibitem{rudin} \textsc{W.~Rudin}, \textsl{Real and Complex Analysis} (McGraw-Hill, Munich, 1987).
\bibitem{shidlovskii} \textsc{A.~B.~Shidlovskii}, \textsl{Transcendental numbers} (Walter de Gruyter, Berlin, 1989).
\bibitem{shimomura} \textsc{S.~Shimomura}, `Growth of the first, the second and the fourth Painlev\'e transcendents',
to appear in Math.~Proc.~Cambridge Philos.~Soc.
\bibitem{steinmetz} \textsc{N.~Steinmetz}, `Value distribution of the
Painlev\'e transcendents', Israel J.~Math.~128 (2002) 29--52.
\bibitem{vojta} \textsc{P.~Vojta}, \textsl{Diophantine Approximations and
Value Distribution Theory}, Lecture Notes in Math~1239
(Springer-Verlag, Berlin, 1987).
\bibitem{weyl} \textsc{H.~Weyl}, \textsl{Meromorphic functions and analytic curves}
(Princeton University Press, Princeton, 1943).
\bibitem{yosida} \textsc{K.~Yosida}, `A generalization of Malmquist's theorem', J.~Math.~9 (1933)
253--256.



\end{thebibliography}
\end{document}